\documentclass[12pt]{article}
\usepackage{CJK}
\usepackage{amsbsy}
\usepackage{amsmath}
\usepackage{amssymb}
\usepackage{epsfig}
\usepackage{multicol}
\usepackage{multirow}
\usepackage{color}
\usepackage{bm}
\usepackage{threeparttable}
\usepackage{float}
\renewcommand{\baselinestretch} {1.3}
\makeatletter 
\def\singlespace{\def\baselinestretch{1}\@normalsize}

\@addtoreset{equation}{section}
\newtheorem{theorem}{Theorem}
\newtheorem{lemma}{Lemma}
\newtheorem{remark}{Remark}

\renewcommand{\hat}{\widehat}

\newcommand{\bg}{\begin{eqnarray}}
\newcommand{\ed}{\end{eqnarray}}
\newcommand{\bgn}{\begin{eqnarray*}}
\newcommand{\edn}{\end{eqnarray*}}

\makeatletter
\def\singlespace{\def\baselinestretch{1}\@normalsize}

\date{\today}
\begin{document}

\title{Strong Asymptotic Properties of Kernel Smoothing Estimation for NA Random Variables with Right Censoring}

\author{
Jian-hua Shi$^{a,b,c,d}$, 
 Jian-sen Xu$^{a}$,  
  Jin-feng Xu$^{a*}$ \\ 
 {\footnotesize $^{a}$School of Mathematics and Statistics,}
          {\footnotesize Minnan Normal University, Zhangzhou, 363000, China }\\
           {\footnotesize $^{b}$ Fujian Key Laboratory of Granular Computing and Applications, Zhangzhou, 363000, China}\\
          {\footnotesize $^{c}$ The Institute of Meteorological Big Data-Digital Fujian, Zhangzhou, 363000, China}\\
            {\footnotesize $^{d}$ Fujian Key Laboratory of Data Science and Statistics, Zhangzhou, 363000, China}    
}

\maketitle

\begin{singlespace}
\begin{footnotetext}
{
Corresponding author: Jin-feng Xu, Minnan Normal University, Zhangzhou, Fujian, China 363000.
E-mail: xjf@mnnu.edu.cn.
}
\end{footnotetext}
\end{singlespace}

\begin{abstract}

 Most studies for negatively associated
 (NA) random variables consider the complete-data situation, which is actually a relatively ideal condition  in practice.
The paper relaxes this condition to  the incomplete-data  setting and  considers kernel smoothing density and hazard function
estimation in the presence of right censoring based on the Kaplan-Meier estimator. We establish the strong asymptotic properties for these two estimators to assess their asymptotic  behavior and justify their practical use.



\end{abstract}

\bigskip
\bigskip

\baselineskip=20pt
\section{Introduction}

\begin{flushleft}

Negatively associated (NA) random sequence is a sequence of dependent random variables, which was first
introduced by Alam and Saxena in 1981 and then delicately studied by Joag-Dev and Proschan in 1983. The definition is  given as follows.

{\textbf{Definition}} \emph{(Joag-Dev and Proschan, 1983) Random sequences \{${{T}_{i}}, 1 <i \leq n$\}  are said to be negatively associated (NA) if for every pair of disjoint subsets ${{B}_{1}}$ and ${{B}_{2}}$ from $\left\{ 1,2,\ldots ,n \right\}$,}
\end{flushleft}
$$\operatorname{cov}\left( {{f}_{1}}({{T}_{i}};i\in {{B}_{1}}),{{f}_{2}}({{T}_{j}};j\in {{B}_{2}}) \right)\le 0,$$
\emph{where there exists the covariance for ${{f}_{1}}(\cdot)$ and ${{f}_{2}}(\cdot)$ with increasing for every variable (or decreasing for every variable). A sequence of random variables $\left\{ {{T}_{i}};i\ge 1 \right\}$ is said to be NA if every
finite subfamily is NA.} 

Clearly,  independent random variable sequences are NA, and  many other non-independent random sequences, for example the random sampling without replacement in a finite population, can also be included in NA category.  Many researchers have studied the properties of NA random variables and obtained many important results. For example, Su et al. (1997) established a probability inequality and some moment inequalities for the partial sum of a NA sequence, which can be used to prove some properties for strictly stationary NA sequences such as weak invariance principle. Shao (2000) proved that most of the well-known  inequalities, such as the Kolmogorov exponential inequality and the Rosenthal
maximal inequality, still hold for NA random variables.
Wu and Chen (2013) presented two strong representation results of the Kaplan-Meier estimator for NA data with censoring, which are most
relevant to the main results in our paper. Zhou and Lin (2015) considered a nonparametric regression model with repeated NA error structures where the wavelet method is used to estimate the regression function.
Thuan and Quang (2016) studied some properties for the NA random variables and
obtained some inequalities including maximal inequality and H$\acute{a}$jek- R$\acute{e}$nyi's type inequality. Tang et al. (2018)  studied the asymptotic normality of the wavelet estimator  in nonparametric regression, where the random errors  are asymptotically NA random variables. Meng (2018) established two general strong laws of large numbers which also involve NA random variables.

Most studies for NA random variable are under a complete-data setting, however, which is actually  a relatively ideal condition  in practice.
In survival analysis, right censoring is often encountered. For detailed discussion of the censoring and its practical relevance, please see Gijbels and Wang (1993), Zhou and Yip(1999), Chen et al. (2015), Qiu et al. (2015), Shi et al. (2018), Ma et al. (2019), Zhang and Zhou (2018) among many others for reference.

Let $(T_i, Y_i), i=1,\ldots ,n, $ denote a sequence of random vectors where $T_i\geq 0$ is the true survival time of interest, which is right censored by the random variable $Y_i\geq 0$. It is assumed that $T_i$ is independent of $Y_i$, but the i.i.d. assumption is not made  for $T_i$'s and $Y_i$'s, which are both NA in our paper. The observations consist of
$({{X}_{i}},{{\delta }_{i}})$, where
\begin{center}
${{X}_{i}}\triangleq \min ({{T}_{i}},{{Y}_{i}})\triangleq{{T}_{i}}\wedge {{Y}_{i}}$~~~  and~~~   ${{\delta }_{i}}\triangleq I({{T}_{i}}\le {{Y}_{i}}), \ i=1,\ldots ,n,$
\end{center}
and  $I(A)$ is the indicator function of the random event $A.$
 For simplicity, assume that $T_i$ have a common continuous marginal distribution function $F(t)\triangleq P({{T}_{i}}\le t)$ and let its survival distribution $S_T(t)\triangleq 1-F(t)$. The random censoring times ${{Y}_{i}}, i=1,\ldots ,n$, being independent of the random variables ${{T}_{i}}$'s , are assumed to have a common distribution function $G(t)\triangleq P({{Y}_{i}}\le t)$  with its survival distribution $S_Y(t)\triangleq 1-G(t)$. Meanwhile, let $L(\cdot)$ be the distribution of the observed variable ${{X}_{i}}$'s, and we write its survival distribution as $\bar{L}(t)\triangleq 1-L(t).$ For any distribution function $H(\cdot)$, we define the left and right endpoints of its support as ${{a}_{H}}$ and ${{\tau }_{H}}$ by  $ {{a}_{H}}\triangleq \inf \{t: H(t)>0\},{{\tau }_{H}}\triangleq \sup \{t: H(t)<1\}$ throughout the  paper.

The distribution function $L(\cdot)$ can be consistently estimated by the empirical distribution function ${{L}_{n}}(t)$, which is defined as follows.
$${{L}_{n}}(t)\triangleq \frac{1}{n}\sum\limits_{k=1}^{n}{I({{X}_{k}} \leq t)}\triangleq 1-\frac{{{Y}_{n}}(t)}{n}\triangleq\frac{{{{\bar{Y}}}_{n}}(t)}{n},$$
where ${{Y}_{n}}(t)\triangleq \sum\limits_{k=1}^{n}{I({{X}_{k}}> t)}$ denotes the number of uncensored and censored observations larger than time $t$, and ${{\bar{Y}}_{n}}(t)\triangleq \sum\limits_{k=1}^{n}{I({{X}_{k}} \leq t)}$.

For drawing nonparametric inference about unknown $F(\cdot)$ based on the censored observations $(X_i,\delta_i),i=1,\ldots ,n,$ we introduce
a stochastic process on $[0,\infty )$ as follows
$${{N}_{n}}(t)\triangleq \sum\limits_{k=1}^{n}{I({{T}_{k}}\le t\wedge {{Y}_{k}})},$$ 
which counts the number of uncensored observations no larger than time $t$. One nonparametric maximum likelihood estimation ${{\hat{F}}_{n}}(\cdot)$ of $F(\cdot)$ is the well-known Kaplan-Meier (K-M) estimator (Kaplan and Meier, 1958), which is commonly used to estimate $F(\cdot)$ for the incomplete data $({{X}_{i}},{{\delta }_{i}})$, i.e.
$$1-{{\hat{F}}_{n}}(x)\triangleq \prod\limits_{t\le x}{(1-\frac{d{{N}_{n}}(t)}{{{Y}_{n}}(t)})},$$
where the jump $d{{N}_{n}}(t)\triangleq {{N}_{n}}(t)-{{N}_{n}}(t-).$

Define the sub-distribution function ${{F}_{*}}(t)\triangleq P({{X}_{1}}\le t,{{\delta }_{1}}=1).$
Since $F(0)=0$,  using integration by parts, we have
$$\begin{aligned}
  & {{F}_{*}}(t)=\int_{0}^{\infty }{\left[ \int_{s\le t \wedge z}{dF(s)} \right]}dG(z)\\ 
 & \quad \quad \,=\int_{0}^{t}{F(z)}dG(z)+\int_{t}^{\infty }{F(t)}dG(z) \\
 & \quad \quad \,=-\int_{0}^{t}{F(z)}dS_Y(z)+F(t)S_Y(t) =\int_{0}^{t}{S_Y(z)}dF(z), \\
\end{aligned}$$
and then
$$d{{F}_{*}}(t)=S_Y(t)dF(t).$$

 It is further assumed that $F(\cdot)$ has a density function $f(\cdot)$. The estimation for the hazard function $h(\cdot)$ is also of substantial interest in survival analysis, which is defined as
$$h(t)\triangleq \frac{d}{dt}\left( -\log S_T(t) \right)=\frac{f(t)}{S_T(t)} \quad for\ S_T(t) > 0.$$
Its correspond cumulative hazard function is defined as
\begin{eqnarray}
H(t)\triangleq \int_0^t {h(s)ds}=\int_{0}^{t}{\frac{d{{F}_{*}}(s)}{\bar{L}(s)}}.  \label{eq1.1}
\end{eqnarray}
The representation (\ref{eq1.1}) of $H(\cdot)$ in terms of ${F}_{*}(\cdot)$ and $\bar{L}(\cdot)$ suggests the empirical estimator for $H(\cdot)$ by
\begin{eqnarray}
{{\hat{H}}_{n}}(t)\triangleq \int_{0}^{t }{\frac{d{{N}_{n}}(s)}{{{Y}_{n}}(s)}}=\int_{0}^{t }{\frac{d{{F}_{*n}}(s)}{{{{\bar{L}}}_{n}}(s)}}, \label{eq1.2}
\end{eqnarray}
where ${{\bar{L}}_{n}}(t)\triangleq 1-{{L}_{n}}(t),$ and
$${{F}_{*n}}(t)\triangleq  \frac{1}{n}\sum\limits_{k=1}^{n}{I({{X}_{k}}\le t,{{\delta }_{k}}=1)}=\frac{{{N}_{n}}(t)}{n}$$
denotes the empirical estimator of ${{F}_{*}(\cdot)}$.

Note that $d{N_n}({X_{(k)}}) = \sum\limits_{j = 1}^n {[{\delta _j}  I({X_j} = {X_{(k)}})]}  \triangleq {\delta _{(k)}},k=1,2,\ldots ,n,$ where ${{X}_{(1)}}\le {{X}_{(2)}}\le \ldots \le {{X}_{(n)}}$ are the order statistics of ${{X}_{1}},{{X}_{2}},\ldots ,{{X}_{n}}$, and ${{\delta }_{(k)}}$ is the concomitant of ${{X}_{(k)}}$.  It can be verified  that the estimators ${{\hat F}_n}(\cdot)$ and ${\hat H_n}(\cdot)$ can be respectively represented as
\begin{eqnarray}
\quad \quad \quad  1 - {{\hat F}_n}(t) = \prod\limits_{{X_{(k)}} \le t} {(1 - \frac{{d{N_n}({X_{(k)}})}}{{n - k + 1}})}  = \prod\limits_{{X_{(k)}} \le t} {(1 - \frac{{{\delta _{(k)}}}}{{n - k + 1}})} ,\label{eq1.3}
\end{eqnarray}
and
\begin{eqnarray}
 {\hat H_n}(t) = \sum\limits_{{X_{(k)}} \le t} {\frac{{{\delta _{(k)}}}}{{n - k + 1}}}. \label{eq1.4}
\end{eqnarray}

In the case of right censoring, the K-M estimator ${{\hat{F}}_{n}}(x)$ and the estimator ${{\hat{H}}_{n}}(x)$ have been generally accepted as a substitute for the usual empirical estimators of distribution function $F(\cdot)$ and the cumulative hazard function $H(\cdot)$, respectively, which help to study other estimators such as the kernel density estimator and the kernel hazard estimator in the following.

A kernel smoothed estimator for $f(\cdot)$ based on ${\hat F}_n(\cdot)$ can be constructed as
\[{\tilde{f}_n}(t)\triangleq b_n^{ - 1}\int_{{a_F}}^{ + \infty } {k(\frac{{t - x}}{{{b_n}}})} d{\hat F_n}(x), \]
where $k(\cdot)$ is a smooth probability kernel function and  $\{ {b_n},n \ge 1\}$ is a sequence of bandwidth tending to zero at appropriate rates.

Similarly, we can also construct a kernel smoothed  estimator for the hazard function $h(\cdot)$ under the NA sampling random variables, which is defined by
\[{\tilde{h}_n}(t) \triangleq b_n^{ - 1}\int_{{a_F}}^{ + \infty } {k(\frac{{t - x}}{{{b_n}}})} d{{\hat H}_n}(x).\]

The estimators ${\tilde{f}_n}(\cdot)$ and ${\tilde{h}_n}(\cdot)$ have attracted the attention of many investigators. For example, Mielniczuk (1986) investigated kernel estimator of a density function using the K-M estimator for censored data. When the data was sampled from $\alpha$-mixing and censoring, Cai (1998) explored the uniform consistency (with rates) and the asymptotic normality of the kernel  estimators  for  density  and  hazard  function. Zhou (1999) successfully established several asymptotic uniformly strong and weak representations for kernel estimators of the density function and the hazard function under left truncation.
Antoniadis et al. (1999) proposed  a wavelet method for estimating density and hazard rate functions from randomly right-censored data.
Some other results, one may refer to Diehl and Stute (1988), Gijbels and Wang (1993), Arcones and Gin$\acute{e}$ (1995), Zhou and Yip (1999),  Lemdani and Ould-Sa\"{i}d (2007), Shen and He (2008) among others.

To present our main results, define
\begin{eqnarray*}
{f_n}(t) &\triangleq& b_n^{ - 1}\int_{{a_F}}^{ + \infty } {k(\frac{{t - x}}{{{b_n}}})} d{F_*}(x), \ \
{h_n}(t) \triangleq b_n^{ - 1}\int_{{a_F}}^{ + \infty } {k(\frac{{t - x}}{{{b_n}}})} dH(x),\\
f_n^*(t) &\triangleq& b_n^{ - 1}\int_{{a_F}}^{ + \infty } {k(\frac{{t - x}}{{{b_n}}})} dF_{*n}(x), \ {F_{*n}}(x) \triangleq \frac{1}{n}\sum\limits_{i = 1}^n {I(T_i\leq x)}.
\end{eqnarray*}

The main purpose of this paper is to study the asymptotic properties of kernel smoothing density estimator $\tilde{f}_n(\cdot)$ and hazard estimator ${\tilde{h}_n}(\cdot)$ based on censoring NA  random variables. Under certain regularity conditions, we establish the  strong asymptotic properties for the two estimators with the convergent rates $O(b_{n}^{-1}{{({{n}^{-1}}\ln n)}^{{1}/{2}\;}})\ a.s.$, where $\{ {b_n},n \ge 1\}$ will be defined in Section 2. Throughout the paper, the sequences of variables  $\{{{T}_{n}};n\ge 1\}$ and $\{{{Y}_{n}};n\ge 1\}$ are all non-negative unless otherwise specified.

\section{Main results and their proofs}
 We first present two lemmas (Wu and Chen, 2013) that will help to prove our theorems.
\begin{lemma} Let $\{{{T}_{n}};n\ge 1\}$ and $\{{{Y}_{n}};n\ge 1\}$ be two sequences of NA random variables. Suppose that the sequences $\{{{T}_{n}};n\ge 1\}$ and $\{{{Y}_{n}};n\ge 1\}$ are independent. Then, for any $0<\tau <{{\tau }_{L}}={{\tau }_{F}}\wedge {{\tau }_{G}},$
\begin{eqnarray}
\underset{0 \leq t \leq \tau }{\mathop{\sup }}\,\left| {{{\hat{F}}}_{n}}(t)-F(t) \right|=O({{({{n}^{-1}}\ln n)}^{{1}/{2}\;}})_{{}}^{{}}~a.s. \label{eq2.2}
\end{eqnarray}
and
\begin{eqnarray}
\underset{0 \leq t \leq \tau }{\mathop{\sup }}\,\left| {{{\hat{H}}}_{n}}(t)-H(t) \right|=O({{({{n}^{-1}}\ln n)}^{{1}/{2}\;}})_{{}}^{{}}~a.s.
\label{eq2.1}
\end{eqnarray}
\end{lemma}

For positive real numbers $x$ and $~t$, write
\[
\eta (x ,t,\delta)\triangleq \int_0^{t \wedge x} {\frac{{d{F_*}(s)}}{{{{\bar L}^2}(s)}}} -\frac{{I(x \le t,\delta  = 1)}}{{\bar L(x)}}.
\]

\begin{lemma} Let $\{{{T}_{n}};n\ge 1\}$ and $\{{{Y}_{n}};n\ge 1\}$ be two sequences of NA random variables. Suppose that the sequences $\{{{T}_{n}};n\ge 1\}$ and $\{{{Y}_{n}};n\ge 1\}$ are independent. Then, for any $0<\tau <{{\tau }_{L}},$
\begin{eqnarray}
{{\hat{F}}_{n}}(t)-F(t)=-\frac{S_T(t)}{n}\sum\limits_{i=1}^{n}{\eta ({{X}_{i}},t,{{\delta }_{i}})}+{{r}_{1n}}(t),
\label{eq2.4}
\end{eqnarray}
and
\begin{eqnarray}
{{\hat{H}}_{n}}(t)-H(t)=-\frac{1}{n}\sum\limits_{i=1}^{n}{\eta ({{X}_{i}},t,{{\delta }_{i}})}+{{r}_{2n}}(t),
\label{eq2.3}
\end{eqnarray}
where $\underset{0 \leq t \leq \tau }{\mathop{\sup }}\,\left| {{r}_{in}}(t) \right|=O({{({{n}^{-1}}\ln n)}^{{1}/{2}\;}})\ a.s., \ i=1,2$.
\end{lemma}

\begin{remark} Note that by the definition of $\eta (x ,t,\delta)$,
\begin{eqnarray*}
   -\frac{1}{n}\sum\limits_{i=1}^{n}{\eta ({{X}_{i}},t,{{\delta }_{i}})}&=& \frac{1}{n}\sum\limits_{i:{{X}_{(i)}}\le t}{\frac{{{N}_{n}}({{X}_{(i)}})-{{N}_{n}}(X_{(i)}-)}{\bar{L}({{X}_{(i)}})}}-\frac{1}{n}\int_{0}^{t}{\frac{\sum\nolimits_{i=1}^{n}{I({{X}_{i}}\ge s)}}{{{{\bar{L}}}^{2}}(s)}}d{{F}_{*}}(s) \\
 &= & \int_{0}^{t}{\frac{1}{\bar{L}(s)}}d{{F}_{*n}}(s)-\int_{0}^{t}{\frac{{{{\bar{L}}}_{n}}(s)}{{{{\bar{L}}}^{2}}(s)}}d{{F}_{*}}(s).
\end{eqnarray*}
Therefore, we can obtain by Lemma 2 that
\begin{eqnarray*}
   {{{\hat{F}}}_{n}}(t)-F(t)
 &  =&(1-F(t))\left[ \int_{0}^{t}{\frac{1}{\bar{L}(s)}}d[{{F}_{*n}}(s)-{{F}_{*}}(s)]-\int_{0}^{t}{\frac{{{{\bar{L}}}_{n}}(s)-\bar{L}(s)}{{{{\bar{L}}}^{2}}(s)}}d{{F}_{*}}(s) \right]\\
 &&+{{r}_{1n}}(t).
\end{eqnarray*}
\end{remark}

One can establish the following lemma by noting that $\{ {X_n};n \ge 1\}$ and $\{ ({X_n},{\delta _n});n \ge 1\}$ are both sequences of NA random variables according to Joag-Dev and Proschan (1983).

\begin{lemma} Under the conditions of Lemma 1, for any $0<\tau <{{\tau }_{L}},$ there are
\begin{eqnarray}
\mathop {\sup }\limits_{0 \leq t \leq \tau  } \left| {{F_{*n}}(t) - {F_*}(t)} \right| = O({({n^{ - 1}}\ln n)^{{1 \mathord{\left/
 {\vphantom {1 2}} \right. \kern-\nulldelimiterspace} 2}}})\ a.s.
 \label{eq2.6}
\end{eqnarray}
and
\begin{eqnarray}
\mathop {\sup }\limits_{0 \leq t \leq \tau  } \left| {{L_n}(t) - L(t)} \right| = O({({n^{ - 1}}\ln n)^{{1 \mathord{\left/
 {\vphantom {1 2}} \right.
 \kern-\nulldelimiterspace} 2}}})\ a.s.
 \label{eq2.5}
\end{eqnarray}
\end{lemma}

\begin{theorem} Under the conditions of Lemma 1, assume that the kernel density $k(t )$ has bounded variation on some finite interval $(r,s)$ with $k(t) = 0$ for $t \notin (r, s)$, where $-\infty< r<0<s<\infty$. Suppose that density distribution $f(\cdot )$ and $g(\cdot )={G}'(\cdot )$ are bounded on the closed interval  $[0,\tau ]$ for some $\tau \in (a_F, {{\tau }_{L}}).$ Then there is
\begin{eqnarray}
   \mathop {\sup }\limits_{0 < t \le \tau } \left| {{\tilde{f}_n}(t) - {{f}_n}(t) - \frac{{f_n^*(t) - Ef_n^*(t)}}{{1 - G(t)}}} \right| = O(b_n^{ - 1}{({n^{ - 1}}\ln n)^{{1 \mathord{\left/
 {\vphantom {1 2}} \right.
 \kern-\nulldelimiterspace} 2}}})\ a.s.,
 \label{eq2.7}
\end{eqnarray}
where the sequence $\{{{b}_{n}};n \ge 1\} $ satisfies $b_n^{ - 1} = o({(n{\ln ^{ - 1}}n)^{{1 \mathord{\left/ {\vphantom {1 2}} \right.
 \kern-\nulldelimiterspace} 2}}})$.
 \end{theorem}

\begin{theorem} Under the conditions of Lemma 1, assume that the kernel density $k(t )$ has bounded variation on some finite interval $(r,s)$ with $k(t) = 0$ for $t \notin (r, s)$, where $-\infty< r<0<s<\infty.$ Suppose that density distribution $f(\cdot )$ and $g(\cdot )={G}'(\cdot )$ are bounded on the closed interval  $[0,\tau ]$ for some $\tau \in (a_F, {{\tau }_{L}}).$ Then there is
\begin{eqnarray}
 \mathop {\sup }\limits_{0 < t \le \tau } \left| {{\tilde{h}_n}(t) - {{ h}_n}(t) - \frac{{f_n^*(t) - Ef_n^*(t)}}{{1 - L(t)}}} \right| = O(b_n^{ - 1}{({n^{ - 1}}\ln n)^{{1 \mathord{\left/
 {\vphantom {1 2}} \right.
 \kern-\nulldelimiterspace} 2}}})\ a.s.,
 \label{eq2.8}
\end{eqnarray}
where the sequence $\{{{b}_{n}};n \ge 1\} $ satisfies $b_n^{ - 1} = o({(n{\ln ^{ - 1}}n)^{{1 \mathord{\left/ {\vphantom {1 2}} \right.
 \kern-\nulldelimiterspace} 2}}})$.
 \end{theorem}

\begin{remark} Theorem 1 and Theorem 2 are key results in studying censored NA sequences, which can be useful in deriving some asymptotic properties for the kernel density estimator $\tilde{f}_n(\cdot)$ and the hazard function estimator $\tilde{h}_n(\cdot)$, respectively. For example, if one can establish the results similar to those in Hall (1981) for NA sequences, then by using Theorem 1, the following proposition will hold.

\textbf{Proposition 1} Suppose  that the sequence $\{{{b}_{n}};n \ge 1\} $ satisfies ${b_n} \to 0$ for $n \to \infty $, and

(a) ${{{{(\ln n)}^2}} \mathord{\left/
 {\vphantom {{{{(\ln n)}^2}} {(n{b_n}\ln \ln n)}}} \right.
 \kern-\nulldelimiterspace} {(n{b_n}\ln \ln n)}} \to 0,$

(b) $n{b_n} \to \infty $  in such a way that
\[\mathop {\lim }\limits_{n \to \infty } \mathop {\sup }\limits_{m:\left| {m - n} \right| \le n\varepsilon } \left| {\frac{{{b_m}}}{{{b_n}}} - 1} \right| \to 0 \ for\ \varepsilon  \to 0,\]
then there will be
\[\mathop {\lim \sup }\limits_{n \to \infty }  \pm {(\frac{{n{b_n}}}{{2\ln \ln n}})^{{1 \mathord{\left/
 {\vphantom {1 2}} \right.
 \kern-\nulldelimiterspace} 2}}}({\tilde{f}_n}(t) - {{ f}_n}(t)) = {[\varphi (f,G)\int {{k^2}(s)ds} ]^{{1 \mathord{\left/
 {\vphantom {1 2}} \right.
 \kern-\nulldelimiterspace} 2}}} \ a.s.\]
 where $\varphi (f,G)$ is some functional for $f(\cdot)$ and $G(\cdot)$.
\end{remark}

For simplicity and without loss of generality, it can be assumed that $a_F=0$ in the following proof.

\emph{Proof of Theorem 1}~~~According to Remark 1, $ {{\tilde{f}}_{n}}(x)-{{{f}}_{n}}(x)$ can be expressed as
\begin{eqnarray}
  {{\tilde{f}}_{n}}(x)-{{{{f}}}_{n}}(x)
 & =&b_{n}^{-1}\int_{0}^{+\infty }{k(\frac{x-t}{{{b}_{n}}})}d\{(1-F(t))[\int_{0}^{t}{\frac{1}{\bar{L}(s)}}d[{{F}_{*n}}(s) -{{F}_{*}}(s)]\nonumber \\
&&-\int_{0}^{t}{\frac{{{{\bar{L}}}_{n}}(s)-\bar{L}(s)}{{{{\bar{L}}}^{2}}(s)}}d{{F}_{*}}(s)]+{{r}_{1n}}(t)\}\nonumber \\
 & = & - b_n^{ - 1}\int_{{0}}^{ + \infty } {\{ (1 - F(t))\int_0^t {\frac{1}{{\bar L(s)}}} d[{F_{*n}}(s) - {F_*}(s)\} } dk(\frac{{x - t}}{{{b_n}}}) \nonumber \\
 && + b_n^{ - 1}\int_{{0}}^{ + \infty } {\{ (1 - F(t))\int_0^t {\frac{{{{\bar L}_n}(s) - \bar L(s)}}{{{{\bar L}^2}(s)}}} d{F_*}(s)\} } dk(\frac{{x - t}}{{{b_n}}}) \nonumber \\
 && - b_n^{ - 1}\int_{{0}}^{ + \infty } {{r_{1n}}(t)} dk(\frac{{x - t}}{{{b_n}}}) \nonumber \\
 & \triangleq & -{{I}_{1}}+{{I}_{2}}-{{I}_{3}}.
  \label{eq3.1}
\end{eqnarray}

Considering ${{I}_{1}}$, we have
\begin{eqnarray*}
  && \int_{0}^{t}{\frac{1}{\bar{L}(s)}}d[{{F}_{*n}}(s)-{{F}_{*}}(s)]
=\frac{{{F}_{*n}}(t)-{{F}_{*}}(t)}{\bar{L}(t)}+\int_{0}^{t}{\frac{{{F}_{*n}}(s)-{{F}_{*}}(s)}{{{{\bar{L}}}^{2}}(s)}}d\bar{L}(s). \end{eqnarray*}

Thus, we have the following formula
\begin{eqnarray}
   {{I}_{1}} & =& b_{n}^{-1}\int_{0}^{+\infty }{(1-F(t))\frac{{{F}_{*n}}(t)-{{F}_{*}}(t)}{\bar{L}(t)}}dk(\frac{x-t}{{{b}_{n}}}) \nonumber \\
 && +b_{n}^{-1}\int_{0}^{+\infty }{(1-F(t))\int_{0}^{t}{\frac{{{F}_{*n}}(s)-{{F}_{*}}(s)}{{{{\bar{L}}}^{2}}(s)}}d\bar{L}(s)}dk(\frac{x-t}{{{b}_{n}}})\nonumber \\
 &\triangleq& {{I}_{11}}+{{I}_{12}}.
  \label{eq3.2}
\end{eqnarray}

Using the partial integration for ${{I}_{11}}$, we get
\begin{eqnarray}
   {{I}_{11}}
 & =& b_{n}^{-1}\int_{0}^{+\infty }{\frac{{{F}_{*n}}(t)-{{F}_{*}}(t)}{1-G(t)}}dk(\frac{x-t}{{{b}_{n}}}) \nonumber \\
 & =& \frac{1}{{{b}_{n}}(1-G(x))}\int_{0}^{+\infty }{[{{F}_{*n}}(t)-{{F}_{*}}(t)]}dk(\frac{x-t}{{{b}_{n}}}) \nonumber \\
 && +\frac{1}{{{b}_{n}}(1-G(x))}\int_{0}^{+\infty }{\frac{{{F}_{*n}}(t)-{{F}_{*}}(t)}{1-G(t)}[G(t)-G(x)]}dk(\frac{x-t}{{{b}_{n}}}) \nonumber \\
 &\triangleq& -\frac{1}{1-G(x)}[f_{n}^{*}(x)-Ef_{n}^{*}(x)]+{{I}_{11}^*}.  \label{eq3.3}
 \end{eqnarray}

Note that the kernel function $k(\cdot)$ is zero outside the interval $(r,s)$ and the fact that $G(\cdot)$ is monotone. Then, when $n$ is large enough, by (2.5) and the definition of $\tau$, applying the change of variable formula,  we have
\begin{eqnarray}
\left| {{I}_{11}^*} \right|
   &\le& \frac{1}{{{b}_{n}}{{(1-G(\tau ))}^{2}}}\underset{0<x\le \tau }{\mathop{\sup }}\,\left| {{F}_{*n}}(t)-{{F}_{*}}(t) \right|\int_{r}^{s}{\left| G(x)-G(x-u{{b}_{n}}) \right|}\left| dk(u) \right| \nonumber \\
 & =& O{{({{({{n}^{-1}}\ln n)}^{{1}/{2}\;}})}_{{}}}~a.s.
 \label{eq3.4}
 \end{eqnarray}
and
 \begin{eqnarray*}
   \left| {{I}_{12}} \right|
 & =&\left| b_{n}^{-1}\int_{r}^{s}{(1-F(x-{{b}_{n}}u))\int_{0}^{x-{{b}_{n}}u}{\frac{{{F}_{*n}}(s)-{{F}_{*}}(s)}{{{{\bar{L}}}^{2}}(s)}}d\bar{L}(s)}dk(u) \right| \\
 & \le& \left| b_{n}^{-1}\int_{r}^{s}{(1-F(x-{{b}_{n}}u))\int_{0}^{x}{\frac{{{F}_{*n}}(s)-{{F}_{*}}(s)}{{{{\bar{L}}}^{2}}(s)}}d\bar{L}(s)}dk(u) \right| \\
 & & +\left| b_{n}^{-1}\int_{r}^{s}{(1-F(x-{{b}_{n}}u))\int_{x}^{x-{{b}_{n}}u}{\frac{{{F}_{*n}}(s)-{{F}_{*}}(s)}{{{{\bar{L}}}^{2}}(s)}}d\bar{L}(s)}dk(u) \right| \\
 & =& b_{n}^{-1}\left| \int_{0}^{x}{\frac{{{F}_{*n}}(s)-{{F}_{*}}(s)}{{{{\bar{L}}}^{2}}(s)}}d\bar{L}(s) \right|\left| \int_{r}^{s}{(1-F(x-{{b}_{n}}u))}dk(u) \right| \\
 && +b_{n}^{-1}\left| \int_{r}^{s}{(1-F(x-{{b}_{n}}u))\int_{x}^{x-{{b}_{n}}u}{\frac{{{F}_{*n}}(s)-{{F}_{*}}(s)}{{{{\bar{L}}}^{2}}(s)}}d\bar{L}(s)}dk(u) \right|.
 \end{eqnarray*}

Again, note that
 \begin{eqnarray*}
  &&  \int_{r}^{s}{(1-F(x-{{b}_{n}}u))dk(u)}
  = -{{b}_{n}}\int_{r}^{s}{k(u)f(x-{{b}_{n}}u)du}.
\end{eqnarray*}

Integrating by parts for ${{I}_{12}}$, we have for $0<\tau <{{\tau }_{L}}$,
 \begin{eqnarray*}
   \left| {{I}_{12}} \right|&\le& b_{n}^{-1}\underset{0<x\le \tau}{\mathop{\sup }}\,\left| {{F}_{*n}}(x)-{{F}_{*}}(x) \right|\left| \int_{0}^{x}{\frac{1}{{{{\bar{L}}}^{2}}(s)}}d\bar{L}(s) \right|{{b}_{n}}\left| \int_{r}^{s}{k(u)f(x-{{b}_{n}}u)}du \right| \\
 && +b_{n}^{-1}\underset{0<x\le \tau}{\mathop{\sup }}\,\left| {{F}_{*n}}(x)-{{F}_{*}}(x) \right| \\
 && \cdot \left| \int_{r}^{s}{(1-F(x-{{b}_{n}}u))\int_{x}^{x-{{b}_{n}}u}{\frac{d[(1-F(s))(1-G(s))]}{{{(1-F(s))}^{2}}{{(1-G(s))}^{2}}}}}dk(u) \right| \\
 & \le& \underset{0<x\le \tau}{\mathop{\sup }}\,\left| {{F}_{*n}}(x)-{{F}_{*}}(x) \right|\left| \left. -\frac{1}{\bar{L}(s)} \right|_{0}^{x} \right|\underset{0<x\le \tau}{\mathop{\sup }}\,f(x)\int_{r}^{s}{\left| k(u) \right|du} \\
 & & +b_{n}^{-1}\underset{0<x\le \tau}{\mathop{\sup }}\,\left| {{F}_{*n}}(x)-{{F}_{*}}(x) \right|\\
 && \cdot \left| \int_{r}^{s}{(1-F(x-{{b}_{n}}u))\int_{x}^{x-{{b}_{n}}u}{\frac{f(s)(1-G(s))+g(s)(1-F(s))}{{{(1-F(s))}^{2}}{{(1-G(s))}^{2}}}ds}}dk(u) \right|.
\end{eqnarray*}

Since density function $f(\cdot )$ and $g(\cdot )$ are bounded in the closed interval $[0,\tau ]$, which means that ${\bar{L}}'(s)=f(s)(1-G(s))+g(s)(1-F(s))$ is also bounded in the interval $[0,\tau ]$, and hence
\begin{eqnarray}
   \left| {{I}_{12}} \right|
 & \le& M\cdot\underset{0<x\le \tau}{\mathop{\sup }}\,\left| {{F}_{*n}}(x)-{{F}_{*}}(x) \right|\left| \frac{1}{\bar{L}(T)}-\frac{1}{\bar{L}(0)} \right|\underset{0<x\le \tau}{\mathop{\sup }}\,f(x) \nonumber \\
 && +b_{n}^{-1}\underset{0<x\le \tau}{\mathop{\sup }}\,\left| {{F}_{*n}}(x)-{{F}_{*}}(x) \right|\int_{r}^{s}{\left| \frac{\underset{0<x\le \tau}{\mathop{\sup }}\,{\bar{L}}'(x)}{{{(1-F(\tau))}^{2}}{{(1-G(\tau))}^{2}}}\int_{x}^{x-{{b}_{n}}u}{ds} \right|}\left| dk(u) \right| \nonumber \\
  & =& M\cdot\underset{0<x\le \tau}{\mathop{\sup }}\,\left| {{F}_{*n}}(x)-{{F}_{*}}(x) \right|\left| \frac{1}{\bar{L}(T)}-\frac{1}{\bar{L}(0)} \right|\underset{0<x\le \tau}{\mathop{\sup }}\,f(x) \nonumber \\
 && +\underset{0<x\le \tau}{\mathop{\sup }}\,\left| {{F}_{*n}}(x)-{{F}_{*}}(x) \right|\left| \frac{\underset{0<x\le \tau}{\mathop{\sup }}\,{\bar{L}}'(x)}{{{(1-F(\tau))}^{2}}{{(1-G(\tau))}^{2}}} \right|\int_{r}^{s}{\left|u \right| \left| dk(u) \right|} \nonumber \\
 & =& O{{({{({{n}^{-1}}\ln n)}^{{1}/{2}\;}})}_{{}}}~a.s.,
  \label{eq3.5}
\end{eqnarray}
where $M$ is some positive constant number.

Thus, combining equations   (\ref{eq3.2}) - (\ref{eq3.5}), we have
\begin{eqnarray}
\mathop {\sup }\limits_{0 < x \le \tau} \left| {{I_1} + \frac{{f_n^*(x) - Ef_n^*(x)}}{{1 - G(x)}}} \right| = O{({({n^{ - 1}}\ln n)^{{1 \mathord{\left/ {\vphantom {1 2}} \right. \kern-\nulldelimiterspace} 2}}})}\ a.s.  \label{eq3.6}
 \end{eqnarray}

On the other hand, similar to the discussion of ${{I}_{12}}$,
\begin{eqnarray*}
   \left| {{I}_{2}} \right|
 & =& \left| {b_n^{ - 1}\int_{{0}}^{ + \infty } {\{ (1 - F(t))[\int_0^x {\frac{{{{\bar L}_n}(s) - \bar L(s)}}{{{{\bar L}^2}(s)}}} d{F_*}(s)} } \right. \\
&& \left. { + \int_x^t {\frac{{{{\bar L}_n}(s) - \bar L(s)}}{{{{\bar L}^2}(s)}}} d{F_*}(s)]\} dk(\frac{{x - t}}{{{b_n}}})} \right| \\
 & \le& b_{n}^{-1}\left| \int_{0}^{x}{\frac{{{{\bar{L}}}_{n}}(s)-\bar{L}(s)}{{{{\bar{L}}}^{2}}(s)}}d{{F}_{*}}(s) \right|\left| \int_{r}^{s}{(1-F(x-{{b}_{n}}u))}dk(u) \right| \\
 &&  +b_{n}^{-1}\left| \int_{r}^{s}{\{(1-F(x-{{b}_{n}}u))\int_{x}^{x-{{b}_{n}}u}{\frac{{{{\bar{L}}}_{n}}(s)-\bar{L}(s)}{{{{\bar{L}}}^{2}}(s)}}d{{F}_{*}}(s)\}}dk(u) \right|\\
 &\triangleq& {{I}_{21}}+{{I}_{22}},
\end{eqnarray*}

where
\begin{eqnarray*}
{{I}_{21}}
 & =& \underset{0\le x\le \tau}{\mathop{\sup }}\,\left| {{{\bar{L}}}_{n}}(x)-\bar{L}(x) \right|\left| \int_{0}^{x}{\frac{f(s)ds}{{{(1-F(s))}^{2}}(1-G(s))}} \right|\left| -{{b}_{n}}\int_{r}^{s}{k(u)f(x-{{b}_{n}}u)du} \right| \\
 & \le& \underset{0\le x\le \tau}{\mathop{\sup }}\,\left| {{{\bar{L}}}_{n}}(x)-\bar{L}(x) \right|\cdot \underset{0\le x\le \tau}{\mathop{\sup }}\,{{f}}(x)\cdot \frac{1}{{{(1-F(\tau))}^{2}}(1-G(\tau))},
\end{eqnarray*}

and
\begin{eqnarray*}
 {{I}_{22}}
 & \le& b_{n}^{-1}\underset{0\le x\le \tau}{\mathop{\sup }}\,\left| {{{\bar{L}}}_{n}}(x)-\bar{L}(x) \right|\\
 &&\cdot \left| \int_{r}^{s}{\{(1-F(x-{{b}_{n}}u))[\int_{x}^{x-{{b}_{n}}u}{\frac{(1-G(s))dF(s)}{{{(1-F(s))}^{2}}{{(1-G(s))}^{2}}}}]\}}\left| dk(u) \right| \right| \\
 & \le& b_{n}^{-1}\underset{0\le x\le \tau}{\mathop{\sup }}\,\left| {{{\bar{L}}}_{n}}(x)-\bar{L}(x) \right|\int_{r}^{s}{\left| \int_{x}^{x-{{b}_{n}}u}{\frac{f(s)ds}{{{(1-F(s))}^{2}}(1-G(s))}} \right|\left| dk(u) \right|} \\
 & \le& \underset{0\le x\le \tau}{\mathop{\sup }}\,\left| {{{\bar{L}}}_{n}}(x)-\bar{L}(x) \right|\int_{r}^{s}{\left| u\right| \left| dk(u) \right|}\underset{0\le x\le \tau}{\mathop{\sup }}\,f(x)\left| \frac{1}{{{(1-F(\tau))}^{2}}(1-G(\tau))} \right|.
\end{eqnarray*}

It can be obtained by (\ref{eq2.6}) that
\begin{eqnarray}
\mathop {\sup }\limits_{0 \le x \le \tau} \left| {{I_2}} \right| = O{({({n^{ - 1}}\ln n)^{{1 \mathord{\left/
 {\vphantom {1 2}} \right.
 \kern-\nulldelimiterspace} 2}}})}\ a.s. \label{eq3.7}
\end{eqnarray}

As for term ${{I}_{3}}$, note that $k(\cdot)$ is of bounded variation, it follows from Lemma 2 that
\begin{eqnarray}
\underset{0<x\le \tau}{\mathop{\sup }}\,\left| {{I}_{3}} \right|=b_{n}^{-1}\underset{0<x\le \tau}{\mathop{\sup }}\,{{r}_{1n}}(t)\int_{0}^{+\infty }{\left| dk(\frac{x-t}{{{b}_{n}}}) \right|}=O(b_{n}^{-1}{{({{n}^{-1}}\ln n)}^{{1}/{2}\;}})_{{}}^{{}}\ a.s.\label{eq3.8}
\end{eqnarray}

This completes the proof by combining (\ref{eq3.1}) and  (\ref{eq3.6})- (\ref{eq3.8}).

\emph{Proof of Theorem 2}
Note by the strong asymptotic expression from (\ref{eq2.3}),
\[{\hat H_n}(t) - H(t) = \int_0^t {\frac{1}{{\bar L(s)}}} d{F_{*n}}(s) - \int_0^t {\frac{{{{\bar L}_n}(s)}}{{{{\bar L}^2}(s)}}} d{F_*}(s) + {r_{2n}}(t),\]
and similarly for the term $I_{11}$, we have
\begin{eqnarray*}
 &&  b_n^{ - 1}\int_0^{ + \infty } {\frac{{{F_{*n}}(t) - {F_*}(t)}}{{1- L(t)}}} dk(\frac{{x - t}}{{{b_n}}}) \\
  &=& \frac{1}{{{b_n}(1- L(x))}}\int_0^{ + \infty } {[{F_{*n}}(t) - {F_*}(t)]} dk(\frac{{x - t}}{{{b_n}}}) \\
 && + \frac{1}{{{b_n}(1- L(x))}}\int_0^{ + \infty } {\frac{{{F_{*n}}(t) - {F_*}(t)}}{{1- L(t)}}[ L(t) - L(x)]} dk(\frac{{x - t}}{{{b_n}}}).
 \end{eqnarray*}

 Then following the proofs of Theorem 1, we can also obtain Theorem 2. This completes our proof.\\
 \\

\small
\baselineskip=8pt

\end{document}